\documentclass[11 pt]{amsart}
\usepackage{epsfig,graphics,amsmath,amssymb,amsthm}

\renewcommand{\contentsname}{Contents\\{\footnotesize\normalfont(A table
of contents should normally not be included)}}
\newtheorem{Theorem}{Theorem}[section]

\newtheorem{Lemma}[Theorem]{Lemma}

\theoremstyle{definition}

\begin{document}

\title[Black holes in the Dynnikov Coordinate Plane]
{Black holes in the Dynnikov Coordinate Plane}

\author{Ferihe Atalan}
\address{Department of Mathematics\\Atilim University\\
06830 Ankara\\ Turkey}
\email{ferihe.atalan@atilim.edu.tr}
\date{\today}
 \pagenumbering{arabic}

\begin{abstract} 
This work gives an application of Dynnikov coordinates in geometric group theory. We present a new approach to the study of subgroups of the mapping class group of a disc with $3$ marked points, $\mathbb{D}_3$, using Dynnikov coordinates.

\vspace{3em} 
\begin{flushright} 
    \footnotesize 
    \textit{Black holes ain't
		as black as
		they are painted.\\}
		Stephen Hawking
\end{flushright}
\end{abstract}

\maketitle
\section{Introduction}

Let $\mathbb{D}_3$ be a disc with $3$ marked points. Let ${\rm Mod}(\mathbb{D}_3)$ denote the mapping class group of $\mathbb{D}_3$, defined as the group of isotopy classes of  orientation-preserving diffeomorphisms that preserve both the boundary and the set of marked points setwise, modulo isotopy through diffeomorphisms satisfying the same conditions. In this case, ${\rm Mod}(\mathbb{D}_3)$ is isomorphic to Artin's braid group $B_3$ modulo its center.

\begin{figure}[h!]
\begin{center}
\hskip-4mm 
    \includegraphics[width=6.2cm]{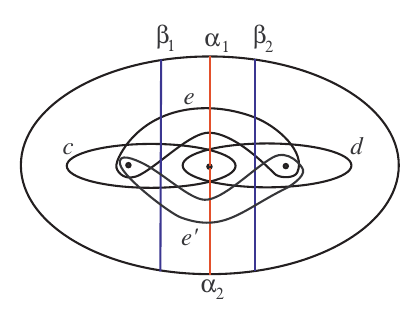}\hskip2mm
    \caption{The curves $c$, $d$, $e$ and $e'$}
    \label{4curves}
  \end{center}
\end{figure}

The dynamics for braid generators on a finitely punctured disc 
in terms of Dynnikov coordinates was studied in \cite{D}, \cite{HY}, \cite{Y}. 
In \cite{AF}, the orbits and dynamics of the action of certain Dehn twists were studied in the Dynnikov coordinate plane. In \cite{AF}, the action of each of the Dehn twists on the Dynnikov coordinate plane $\mathbb{Z}^2$  is piecewise-linear and area-preserving (see \cite{AF}). These actions have an interesting geometric meaning as piecewise linear $\mathbb{Z}^2$-automorphisms preserving the shape of the linearity border fan. This geometric property provides a new approach to prove a classical result: the group generated by Dehn twists $t_i$ for $i=1, \ldots, k$ is isomorphic to the free group $F_k$ of rank $k$ under certain conditions. Such studies can be found in \cite{H, I, HT, FM, K} and \cite{DMYA}. In \cite{DMYA}, it is shown that in $\mathbb{D}_3$ two Dehn twists $t_c$ and $t_d$ generate $F_2$ using Dynnikov matrices (see \cite{Y} for the definition of Dynnikov matrices), where $c$ and $d$ are simple closed curves as shown in Figure\,\ref{4curves}. In this note, using the geometric property derived from the results in \cite{AF}, we provide a geometric proof of the following result, which, to the best of our knowledge, appears to be new.

\begin{Theorem}\label{Thm}
Let $c$, $d$, $e$, and $e'$ be simple closed curves as in depicted Figure\,\ref{4curves}. Then the group generated by the squares of Dehn twists $t^2_c$, $t^2_d$, $t^2_e$, and $t^2_{e'}$  is isomorphic to the free group of rank\,4.
\end{Theorem}

\section{Preliminaries}\label{Prelim}

Let $x$ be a simple closed curve in $\mathbb{D}_3$. If $x$ bounds a disc containing exactly two marked points, we call it an essential simple closed curve. A multi-curve $u$ is a nonempty finite union of mutually disjoint essential simple closed curves in $\mathbb{D}_3$. 
We denote the set of isotopy classes of essential simple closed curves and multi-curves by $\mathcal{C}$ and  $\mathcal{U}$, respectively. We note that any two essential simple closed curves in $\mathbb{D}_3$ intersect each other at least two points.

Throughout the paper, the same notation will be used for curves and their isotopy classes.

For a multi-curve $u$, there is a bijective map $\rho:\mathcal{U} \to \mathbb{Z}^{2} \setminus \{0\}$ defined by $\rho(u)=(a,b)$, where $$a=\frac{\alpha_{2}(u) - \alpha_{1}(u)}{2}  \,\,\, \, \,\,\,  b=\frac{\beta_1(u) - \beta_{2}(u)}{2}.$$
Here, $\alpha_{i}(u)$ and  $\beta_i(u)$ denote the minimum intersection number of $u$ with the arcs $\alpha_{i}$ and  $\beta_i$, respectively, for $i=1,2$ (see \cite{D, HY}).

We observe that the Dynnikov coordinates of $c$, $d$, $e$ and $e'$ in Figure\,\ref{4curves} are $(0,1)$, $(0,-1)$, $(-1, 0)$ and $(1,0)$, respectively.

In \cite{AF}, the action of the pair of Dehn twists $t_c, t_d$ and $t_e, t_{e'}$ on the Dynnikov coordinate plane is given in Theorem\,3.1 and Theorem\,3.6, and also their orbits are shown in Figure\,5 and Figure\,8, respectively.
From Theorem\,3.1, we conclude that positive iteration of $t_c$ (resp. $t_d$) leads to sequences that flow into the second (resp. fourth) quadrant. These quadrants are invariant under the respective Dehn twist. Let us call these quadrants \textit{black holes}.

In this context, let us recall that Yurtta\c{s} also introduced Dynnikov regions to study the dynamics of pseudo-Anosov braids in \cite{Y}. She noted that similar
regions were considered for any orientable surface in \cite{M} using train tracks.

\subsection{Linearity regions} 
\begin{figure}[h!]
\begin{center}
    \includegraphics[width=9.2cm]{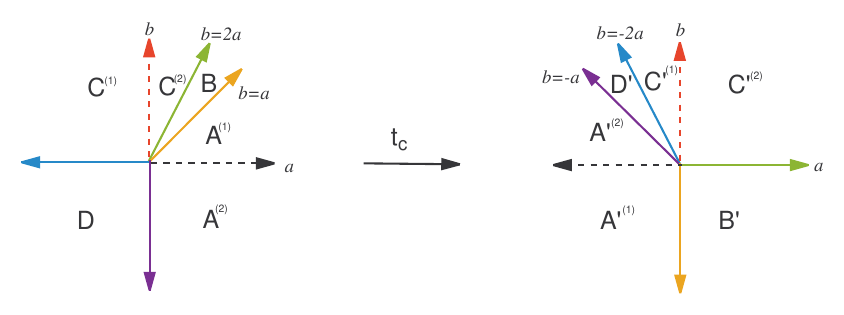}\hskip2mm
   \includegraphics[width=9.2cm]{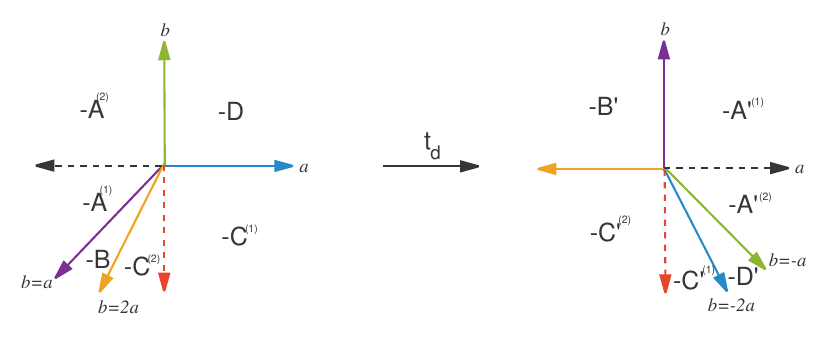}
    \caption{The linearity regions for the $t_c$ and $t_d$-action}
    \label{Fan1-2}
  \end{center}
\end{figure}

\begin{figure}[h!]
 \begin{center}
\includegraphics[width=9.2cm]{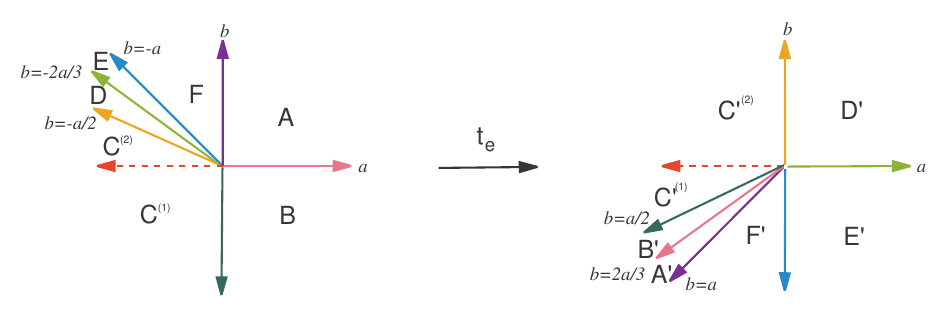}\hskip3mm    
\includegraphics[width=9.2cm]{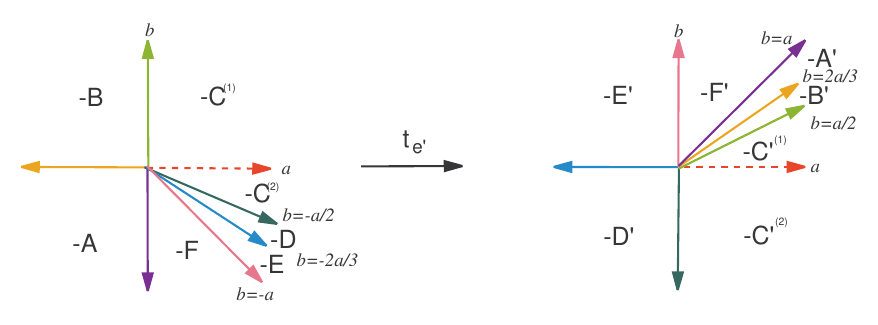}
    \caption{Linearity regions for the action of $t_e$ and $t_{e'}$.} 
		\label{fan3-4}
  \end{center}
\end{figure}

Figure\,\ref{Fan1-2} and Figure\,\ref{fan3-4} show the  linearity regions for the action of Dehn twists $t_c, t_d$ and $t_e, t_{e'}$, and the images of these regions. Namely, the image of region
$B$ is denoted by $B'$, etc.

 \vspace{1cm}

\section{Proof of Theorem\,\ref{Thm}}

\begin{Lemma}(Ping pong lemma \cite{FM})\label{PP-lemma}
Let $G$ be a group acting on a set $\mathcal{C}$. Let $g_1, \ldots, g_k$  be elements of $G$. Suppose that there are nonempty, disjoint subsets $I_1, \ldots, I_k$ of $\mathcal{C}$ with the property that, for each $i$ and each $j \neq i$, we have $g_i^{p}(I_j) \subset I_i$ for every integer $p \neq 0$. Then the group generated by $g_1, \ldots, g_k$ is a free group $F_k$ of rank $k$.
\end{Lemma} 

\subsection{Notation}
Let $(a_l, b_l)$ denote the Dynnikov coordinates of $l \in \mathcal{C}$. 

Let us define sets $I_1$, $I_2$, $I_3$, and $I_4$ as follows: \\
$$I_1 = \{\gamma \in \mathcal{C} \,|\,  b_\gamma \geq  |a_\gamma| \, \}$$  
$$I_2 = \{\gamma \in \mathcal{C} \,|\,  b_\gamma \leq -|a_\gamma| \,\}$$ 
$$I_3 = \{\gamma \in \mathcal{C} \,|\, a_\gamma < 0,\,\, a_\gamma < b_\gamma < -a_\gamma \}$$ 
$$I_4 = \{\gamma \in \mathcal{C} \,|\, a_\gamma > 0,\,\, -a_\gamma < b_\gamma < a_\gamma \}$$

\begin{figure}[hbt]
\hskip-30mm  \includegraphics[width=10cm]{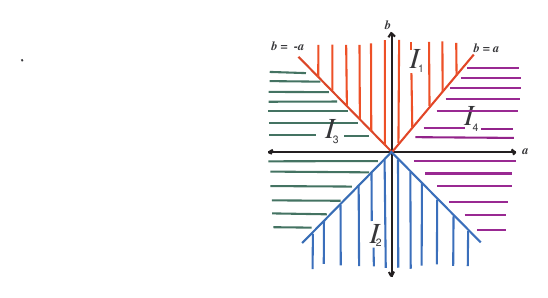}
    \caption{The sets $I_1$,  $I_2$, $I_3$, $I_4$}
    \label{F4}
\end{figure}

It is obvious that these sets are disjoint. They are also nonempty; because, the Dynnikov coordinates of the curves $c$, $d$, $e$ and $e'$ in the Figure\,\ref{4curves} are $(0,1)$, $(0,-1)$, $(-1,0),$ and $(1,0),$ respectively. Since $a_c=0,\,\,b_c = 1 > 0$, so, $c \in I_1$ and  $a_d =0,\,\, b_d = -1 < 0$, so, $d \in I_2$. Similarly, $e \in I_3$ and $e' \in I_4$.  

By applying Lemma\,\ref{PP-lemma}, the proof reduces to verifying that $t_c^{2p}
(I_2 \cup I_3 \cup I_4) \subset I_1$, 
$t_d^{2p}(I_1 \cup I_3 \cup I_4) \subset I_2$, 
$t_e^{2p} (I_1 \cup I_2 \cup I_4) \subset I_3$ and 
$t_{e'}^{2p} (I_1 \cup I_2 \cup I_3) \subset I_4$ for $p \neq 0$. Because of symmetry, we can assume that $p>0$. Moreover, again, because of symmetry, it is sufficient to show that 
$$t_c^{2p} (I_2 \cup I_3 \cup I_4) \subset I_1 and$$ 
$$t_e^{2p} (I_1 \cup I_2 \cup I_4) \subset I_3.$$ Additionally, we use \textit{the black holes} that are invariant under the associated Dehn twists: The second quadrant is invariant under the Dehn twist $t^{p}_c$ for $p>0$ as $b$-coordinates increase, whereas the fourth quadrant is invariant under the Dehn twist $t^{p}_d$ for $p>0$, as $b$-coordinates decrease. On the other hand, the third quadrant is invariant under the Dehn twist $t^{p}_e$ for $p>0$ as $a$- coordinates decrease, whereas the first quadrant is invariant under the Dehn twist $t^{p}_{e'}$ for $p>0$ as $a$-coordinates increase. For $p<0$, consider symmetry with respect to the associated axes.

\begin{figure}[hbt]
  \begin{center}
  \includegraphics[width=13.2cm]{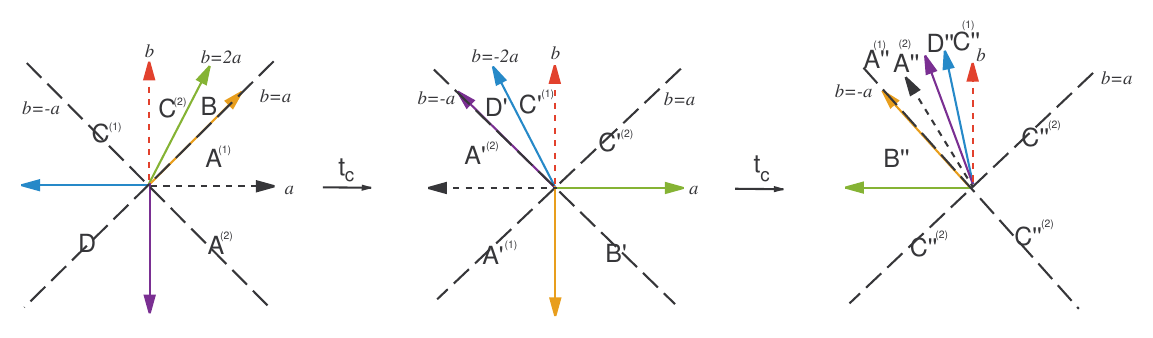}
    \caption{The proof of $t_c^{2}
(I_2 \cup I_3 \cup I_4) \subset I_1$ 
    , i.e., $C''^{(1)} \cup D''  \cup A''^{(2)} \cup A''^{(1)} \subset I_1$} 
		\label{prooft_c}
  \end{center}
\end{figure}

\begin{figure}[hbt]
  \begin{center}
  \includegraphics[width=13.2cm]{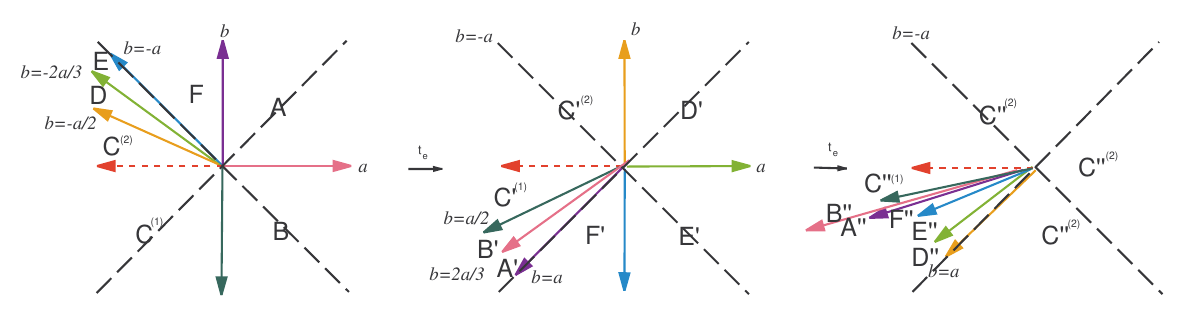}
    \caption{The proof of $t_e^{2} (I_1 \cup I_2 \cup I_4) \subset I_3$, i.e., $C''^{(1)} \cup B'' \cup A'' \cup F'' \subset I_3$ }
    \label{prooft_e}
  \end{center}
\end{figure}

The proof for the positive powers is complete due to the invariance of the black holes. The proof for the negative powers is similar.

\qed

\textbf{Acknowledgements.} I would like to thank the Max Planck Institute for
Mathematics in Bonn for its hospitality, excellent working conditions, and support in summer 2024, when this work was initiated. I would like to thank S. Finashin for useful remarks. I would also like to thank \"{O}. Yurtta\c{s} for some references. Finally, the author is grateful to her mother, Z. Atalan, for her support in summer 2024.


\bigskip
\providecommand{\bysame}{\leavevmode\hboxto3em{\hrulefill}\thinspace}


\begin{thebibliography}{1}

\bibitem{AF} F. Atalan and S. Finashin, \textit{Recognising conjugacy classes of Dehn twists on $\mathbb{D}_3$}, arXiv:2603.14582.

\bibitem{DMYA} E. Dalyan, E. Medeto\~{g}ullari, S.\"{O}. Yurtta\c{s} and F. Atalan, \textit{Generating the free group of rank two with Dynnikov 
Coordinates}, Cumhuriyet Sci. J., 47(2) (2026) 356–360. 

\bibitem{D} I.A. Dynnikov, \textit{On a Yang-Baxter mapping and the Dehornoy ordering}, Russian Mathematical  Surveys,  57(3) (2002) 592-594.

\bibitem{FM} B. Farb and D. Margalit, \textit{A primer on mapping class groups},
Princeton University Press, New Jersey, 2012.

\bibitem{HY} T. Hall and S.\"{O}. Yurtta\c{s}, \textit{On the topological entropy of families of braids}, Topology and its Applications, 156 (2009) 1554–1564.

\bibitem{HT} H. Hamidi-Tehrani, \textit{Groups generated by positive multi-twists and the fake lantern problem}, Algebraic and Geometric Topology, 2(2) 
(2002) 1155-1178.

\bibitem{H} S. P. Humphries, \textit{Free products in mapping class groups generated by Dehn twists} Glasgow Mathematical Journal, 31(2) (1989) 213-218.

\bibitem{I}	A. Ishida,  \textit{The structure of subgroup of mapping class groups generated by two Dehn twists}, Proceedings of the Japan Academy, Series A, Mathematical Sciences, 72(10) (1996) 240–241.

\bibitem{K}	S. Kolay, \textit{Subgroups of the mapping class group of the torus generated by powers of Dehn twists}, arXiv:1909.07360 (2019).

\bibitem{M} D. Margalit, B. Strenner, S. J. Taylor, and S.\"{O}. Yurtta\c{s}, \textit{Quadratic-time computations for pseudo-anosov mapping classes}, Forum of Mathematics, Sigma, 14:e114 (2026), 1-44 

\bibitem{Y1} S.\"{O}. Yurtta\c{s}, \textit{Dynnikov coordinates and pseudo -Anosov braids}, Ph.D. thesis, University of Liverpool, 2011.

\bibitem{Y} S.\"{O}. Yurtta\c{s}, \textit{Dynnikov and train track transition matrices of pseudo-anosov braids}, Discrete and Continuous Dynamical Systems 
36(1) (2016) 541-570.


\end{thebibliography}
\end{document}